# A Relaxation Theorem for Differential Inclusions with Applications to Stability Properties

Brian P. Ingalls*

Control and Dynamical Systems, California Institute of Technology, CA

ingalls@cds.caltech.edu

Eduardo D. Sontag†

Dept. of Mathematics, Rutgers University, NJ

sontag@math.rutgers.edu

Yuan Wang‡

Dept. of Mathematics, Florida Atlantic University, FL

ywang@math.fau.edu

**Abstract**

The fundamental Filippov–Ważwski Relaxation Theorem states that the solution set of an initial value problem for a locally Lipschitz inclusion is dense in the solution set of the same initial value problem for the corresponding relaxation inclusion on compact intervals. In our recent work, a complementary result was provided for inclusions with finite dimensional state spaces which says that the approximation can be carried out over non-compact or infinite intervals provided one does not insist on the same initial values. This note extends the infinite-time relaxation theorem to the inclusions whose state spaces are Banach spaces. To illustrate the motivations for studying such approximation results, we briefly discuss a quick application of the result to output stability and uniform output stability properties.

---

*Research partially carried out while the author was at Rutgers University, supported in part by US Air Force Grant F49620-98-1-0242

†Supported in part by US Air Force Grant F49620-01-1-0063

‡Supported in part by NSF Grant DMS-0072620. Research partially carried out while the author visited the Institute of Systems Science, Chinese Academy of Sciences.

# 1 Introduction

In the study of control systems, it is very often that a system can be modeled by differential inclusions, especially when a system is subject to external signals such as disturbances. In this work we study the approximation of a differential inclusion of the type:

$$\dot{x}(t) \in F(t, x(t)), \tag{1.1}$$

where the trajectories $x(\cdot)$ take values in a separable Banach space $X$, and the set-valued function $F$ is a locally Lipschitz map with nonempty compact values. For such systems of differential inclusions, it is desirable to have $F$ convex and closed. One of the main reasons is that when $F$ is closed and convex, the solutions of the differential inclusion form a closed set in some topology. When $F$ does not possess the convexity property, it may be helpful to consider the corresponding relaxed inclusion defined by

$$\dot{x}(t) \in \operatorname{clco} F(t, x(t)), \tag{1.2}$$

where clco stands for closed convex hull. It is thus interesting to understand how well a trajectory of the relaxation of (1.2) can be approximated by trajectories of the inclusion (1.1). The fundamental result in this aspect is the Filippov–Ważewski Theorem (cf. [2, 3, 7, 8, 9]), which states that the solution set of (1.1) is dense in the solution set of (1.2) in the uniform topology on compact intervals. In particular, the Filippov–Ważewski Theorem says that every trajectory of (1.2) defined on a finite interval can be approximated by trajectories of (1.1) with the same initial condition as the given trajectory of (1.2).

In the recent work [5], a complementary result was presented in the special case when $X$ is of finite dimension. Roughly speaking, it was shown that in the finite dimensional case, the solution set of initial value problems of the type

$$\dot{x}(t) \in F(t, x(t)), \qquad x(0) = \xi_1, \tag{1.3}$$

is dense in the solution set of initial value problems of the type

$$\dot{x}(t) \in \operatorname{clco} F(t, x(t)), \qquad x(0) = \xi_2, \tag{1.4}$$

in the "$C^0$ Whitney topology" on the *infinite* interval $[0, \infty)$. More specifically, the result in [5] provides the existence of trajectories which are approximations in weighted norms on $[0, \infty)$, for example $|f| := \sup_{t \geq 0}\{|f(t)| e^t\}$. Indeed, given any $r : \mathbb{R}_{\geq 0} \to \mathbb{R}_{>0}$, there is an approximation in the norm $|f| := \sup_{t \geq 0}\{|f(t)| r(t)\}$. This is achieved by demanding that the approximation lie in a tube around the reference trajectory which has possibly vanishing radius. In this note, we generalize this result to the case when $X$ is infinite dimensional.

We remark that the approximation result in this note is not a generalization of the Filippov–Ważewski theorem. Given a trajectory of the relaxation (1.4), our result does not guarantee the existence of an approximating trajectory of (1.3) with $\xi_1 = \xi_2$, but rather only with $\xi$ arbitrarily close to $\xi_2$. In fact, it was shown in [5] by a counterexample that one cannot achieve an approximation on the infinite interval with the same initial condition when the error function $r(t)$ is required to decay to 0 as $t \to \infty$. In this note, the example is modified to show that this is still the case even if the error function $r(t)$ is a constant.

Our motivation for this work originated from the study of stability and uniform stability properties of systems with disturbances. Our efforts in this area started in [10] where it was shown that a differential inclusion $\dot{x} \in F(x)$ is globally asymptotically stable if and only if it is uniformly globally asymptotically stable, provided that the state space $X$ is finite dimensional and the set-valued map $F$ admits a parameterization of the form $F(x) = \{f(x,u) : u \in U\}$ where $f(\cdot, \cdot)$ is locally Lipschitz and $U$ is compact. (See [1] for a more general result). As a quick application of the approximation result in this note, we show the result on stability can be extended to the case of output stability.

## 2 Preliminaries

For each $T > 0$, let $\mathcal{L}[0,T]$ be the $\sigma$-field of Lebesgue measurable subsets of $[0,T]$. Let $X$ be a separable Banach space, whose norm is denoted simply by $|\cdot|$. Let $\mathcal{P}(X)$ denote the family of all nonempty closed subsets of $X$. We use $\mathcal{B}(X)$ for the family of Borel subsets of $X$.

For each interval $\mathcal{I} \subseteq [0, \infty)$, let $L^1(\mathcal{I}, X)$ be the Banach space of Bochner integrable functions $u : \mathcal{I} \to X$ with norm $\|u\| = \int_{\mathcal{I}} |u(t)| \, dt$, and let $L^1_{\text{loc}}(\mathcal{I}, X)$ be the corresponding space of locally integrable functions. Let $AC(\mathcal{I}, X)$ be the Banach space of absolutely continuous functions $u : \mathcal{I} \to X$ with the norm $\|u\|_{AC} = |u(0)| + \|\dot{u}\|$.

We define the distance from a point $\xi \in X$ to a set $K \in \mathcal{P}(X)$ as

$$d(\xi, K) := \inf\{|\xi - \eta| : \eta \in K\}.$$

For a set $A \in \mathcal{P}(X)$, let $B(A, r)$ denote the set $\{\xi \in X : d(\xi, A) \leq r\}$. For singleton $A = \{\xi\}$ we write $B(\xi, r)$. For each set $A$ and each constant $c \in \mathbb{R}$, we denote $cA = \{c\xi : \xi \in A\}$.

**Definition 2.1** The *Hausdorff distance* between two sets $K, L \in \mathcal{P}(X)$ is defined as

$$d_H(K, L) := \max\left\{\sup_{\xi \in K} d(\xi, L), \sup_{\eta \in L} d(\eta, K)\right\}.$$

**Definition 2.2** Let $\mathcal{O}$ be an open subset of $X$. Let $\mathcal{I} \subseteq \mathbb{R}_{\geq 0}$ be an interval. The set-valued map $F : \mathcal{I} \times X \to \mathcal{P}(X)$ is said to be *locally Lipschitz* on $\mathcal{O}$ if, for each $\xi \in \mathcal{O}$, there exists a neighbourhood $U \subset \mathcal{O}$ of $\xi$ and a $k_U \in L^1(\mathcal{I}, \mathbb{R})$ so that for any $\eta, \zeta$ in $U$,

$$d_H(F(t, \eta), F(t, \zeta)) \leq k_U(t) |\eta - \zeta| \qquad \text{a.e. } t \in \mathcal{I}.$$

**Definition 2.3** Let $\mathcal{I} \subseteq [0, \infty)$ be an interval. A function $x : \mathcal{I} \to X$ is said to be a *solution of the differential inclusion*

$$\dot{x}(t) \in F(t, x(t)) \tag{2.5}$$

if it is absolutely continuous and satisfies (2.5) for almost every $t \in \mathcal{I}$.

For $T > 0$, a function $x : [0, T) \to X$ is called a *maximal solution of the differential inclusion* if it does not have an extension which is a solution in $X$. That is, either $T = \infty$ or there does not exist a solution $y : [0, T_+] \to X$ with $T_+ > T$ so that $y(t) = x(t)$ for all $t \in [0, T)$. □

The next result on continuous selections of approximations of a relaxed trajectory over finite intervals was derived in [5] based on the Filippov-Ważewski relaxation theorem (see e.g., [9]).

**Lemma 2.4** Let $T > 0$, $\xi_0 \in X$, and a set-valued map $F : [0, T] \times X \to \mathcal{P}(X)$ be given. Consider, for $t \in [0, T]$, the solutions of the differential inclusion

$$\dot{x} \in F(t, x), \tag{2.6}$$

and the solutions of the initial value problem

$$\dot{x} \in \text{ clco } F(t, x), \quad x(0) = \xi_0. \tag{2.7}$$

Suppose $z : [0, T] \to X$ is a solution of (2.7), and let $\varepsilon > 0$ be given. Let

$$\mathcal{T} := \{\xi \in X \ : \ |\xi - z(t)| \leq \varepsilon \text{ for some } t \in [0, T]\},$$

the $\varepsilon$-tube around the image of $z$. Then, if $F$ satisfies

(H1) $F$ is $\mathcal{L}[0, T] \otimes \mathcal{B}(X)$ measurable;

(H2) there exists $k_0 \in L^1([0, T], \mathbb{R})$ such that for any $\xi, \eta \in B(\mathcal{T}, 1)$

$$d_H(F(t, \xi), F(t, \eta)) \leq k_0(t) |\xi - \eta| \quad \text{a.e. } t \in [0, T];$$

(H3) there exists $\alpha \in L^1([0, T], \mathbb{R})$ such that for each $\xi \in B(\mathcal{T}, 1)$

$$\sup\{|\zeta| \ : \ \zeta \in F(t, \xi)\} \leq \alpha(t) \quad \text{a.e. } t \in [0, T];$$

it follows that there exists a $\delta > 0$ and a function $x : [0, T] \times V \to X$, where $V := B(\xi_0, \delta)$ such that

(a) For every $\eta \in V$, the function $t \mapsto x(t, \eta)$ is a solution of the initial value problem

$$\dot{x} \in F(t, x), \quad x(0) = \eta, \quad \text{for } t \in [0, T]; \tag{2.8}$$

(b) the map $\eta \mapsto x(\cdot, \eta)$ is continuous from $V$ into $AC([0, T], X)$;

(c) for each $\eta \in V$,

$$|z(t) - x(t, \eta)| < \varepsilon$$

for all $t \in [0, T]$. □

# 3 Approximations of Relaxed Trajectories

In this section we present the main result of this work.

**Theorem 3.1** *Let $0 < T \leq \infty$, and suppose $F : [0, T) \times X \to \mathcal{P}(X)$ satisfies the following properties.*

*(H1′) $F$ is $\mathcal{L}[0, T) \otimes \mathcal{B}(X)$ measurable;*

*(H2′) for each $R > 0$, there exists $k_R \in L^1_{\text{loc}}([0, T), \mathbb{R})$ such that for any $\xi$, $\eta \in B(0, R)$*

$$d_H(F(t, \xi), F(t, \eta)) \leq k_R(t) |\xi - \eta| \qquad \text{a.e. } t \in [0, T);$$

*(H3′) for each $R > 0$, there exists $\alpha_R \in L^1_{\text{loc}}([0, T), \mathbb{R})$ such that for each $\xi \in B(0, R)$*

$$\sup\{|\zeta| : \zeta \in F(t, \xi)\} \leq \alpha_R(t) \qquad \text{a.e. } t \in [0, T);$$

*Fix $\xi \in X$ and let $z : [0, T) \to X$ be a solution of*

$$\dot{x} \in \text{ clco } F(t, x), \qquad x(0) = \xi.$$

*Let $r : [0, T) \to \mathbb{R}_{>0}$ be a continuous function. Then, there exists some $\eta^0 \in B(\xi, r(0))$ and a solution $\gamma : [0, T) \to X$ of*

$$\dot{x} \in F(t, x), \quad x(0) = \eta^0,$$

*which satisfies*

$$|z(t) - \gamma(t)| \leq r(t) \qquad \forall\, t \in [0, T).$$

*Proof.* Let $\{T_k\}$ be any strictly increasing sequence of times so that $T_0 = 0$ and $T_k \to T$ as $k \to \infty$, and define $\{r_k\}$ by

$$r_k = \min\{r(s) : s \in [T_k, T_{k+1}]\}.$$

Without loss of generality, we assume that $\{r_k\}$ is nonincreasing. Let $z_k = z(T_k)$. On each interval $[T_{k-1}, T_k]$, we will consider the differential inclusions in backward time:

$$\dot{x} \in -F(T_k - t, x), \quad t \in [0, T_k - T_{k-1}], \tag{3.9}$$

and

$$\dot{x} \in \text{clco}(-F(T_k - t, x)), \quad x(0) = z_k, \quad t \in [0, T_k - T_{k-1}] \tag{3.10}$$

Let $\delta_0 = r_0$. We will prove by induction the following: for each $k \geq 1$, there exist $0 < \delta_k \leq r_k$ and a map $x_k : [0, T_k - T_{k-1}] \times V_k$, where $V_k = B(z_k, \delta_k)$, such that the following holds:

- for every $\eta \in V_k$, the function $t \to x_k(t, \eta)$ is a solution of

$$\dot{x} \in -F(T_k - t, x), \qquad x(0) = \eta$$

for $t \in [0, T_k - T_{k-1}]$;

- the map $\eta \mapsto x_k(\cdot, \eta)$ is continuous from $V_k$ into $AC([0, T_k - T_{k-1}], X)$, and in particular, the map $\varphi_k : \eta \mapsto x_k(T_k - T_{k-1}, \eta)$ is continuous on $V_k$;

- for each $\eta \in V_k$,
$$|x_k(t, \eta) - z(T_k - t)| \leq \delta_{k-1} \qquad \forall\, t \in [0, T_k - T_{k-1}],$$
and in particular, $x_k(T_k - T_{k-1}, \eta) \in B(z_{k-1}, \delta_{k-1})$; and

- for each $1 \leq i \leq k$, define inductively $\zeta_k^{k-i} = \varphi_{k-i+1}(\zeta_k^{k-i+1})$ with $\zeta_k^k = z_k$. Then
$$\left|\zeta_{k+1}^j - \zeta_k^j\right| \leq \frac{\delta_j}{2^k} \qquad \forall\, j \leq k. \tag{3.11}$$

Consider the case when $k = 1$. Note that $z(T_1 - t)$ is a solution of (3.10) with $k = 1$. Further, the hypotheses of Lemma 2.4 are satisfied by the map $F$, and hence, $-F$, since there exists $R$ large enough so that $B(0, R)$ contains the image of $z(T_1 - t)$ over $t \in [0, T_1]$ (so that the local conditions assumed in the hypotheses can be made global by redefining the function outside of such a ball). Applying the Lemma with $\varepsilon = \delta_0$, it follows that there exists a $0 < \delta_1 < r_1$ and a function $x_1 : [0, T_1] \times V_1 \to AC([0, T_1], X)$, where $V_1 := B(z(T_1), \delta_1)$, which satisfies

- For every $\eta \in V_1$, the function $t \mapsto x_1(t, \eta)$ is a solution of
$$\dot{x} \in -F(T_1 - t, x), \qquad x(0) = \eta,$$
for $t \in [0, T_1]$;

- the map $\eta \mapsto x_1(\cdot, \eta)$ is continuous from $V_1$ into $AC([0, T_1], X)$;

- for each $\eta \in V_1$,
$$|x_1(t, \eta) - z(T_1 - t)| < \delta_0 \qquad \forall\, t \in [0, T_1];$$
and consequently,

- with $\zeta_0^0 = z_0$, $\zeta_1^0 = x_1(T_1, z_1)$, it holds that $|\zeta_1^0 - \zeta_0^0| < \delta_0$.

Let $k \geq 1$. Suppose for each $1 \leq j \leq k$, there exist $0 < \delta_j < r_j$ and
$$x_j : [0, T_j - T_{j-1}] \times V_j \to X,$$
where $V_j = B(z_j, \delta_j)$, such that

- for every $\eta \in V_j$, the function $t \mapsto x_j(t, \eta)$ is a solution of
$$\dot{x} \in -F(T_j - t, x), \quad x(0) = \eta,$$
for $t \in [0, T_j - T_{j-1}]$;

- the map $\eta \mapsto x_j(\cdot, \eta)$ is continuous into $AC([0, T_j - T_{j-1}], X)$, and so the map $\varphi_j : V_j \to X$ defined by $\varphi_j(\eta) = x_j(T_j - T_{j-1}, \eta)$ is continuous on $V_j$;

- for each $\eta \in V_j$,
$$|z(T_j - t) - x_j(t, \eta)| < \delta_{j-1} \qquad \forall t \in [0, T_{j+1} - t_j];$$

and

- for the sequence $\{\zeta_i^j\}_{i \geq j}$ defined by $\zeta_j^j = z_j$, $\zeta_i^j = \varphi_{j+1}(\zeta_i^{j+1})$, it holds that
$$\left|\zeta_{i+1}^j - \zeta_i^j\right| < \frac{\delta_j}{2^i}$$
for all $0 \leq j \leq i \leq k - 1$.

Below we produce $\delta_{k+1}$, $x_{k+1}$, and $\{\zeta_{k+1}^j\}$ for $0 \leq j \leq k+1$.

Since, for each $j = 1, \ldots, k$, $\varphi_j$ is continuous on $V_j$, and $\varphi_j(V_j) \subseteq (V_{j-1})$, there exists some $0 < \widehat{\delta}_k < \frac{\delta_k}{2^k}$ such that for any $\eta \in B(z_k, \widehat{\delta}_k)$ and any $0 \leq j \leq k - 1$, it holds that
$$\begin{aligned}
&\left|(\varphi_{j+1} \circ \cdots \circ \varphi_k)(\eta) - \zeta_k^j\right| \\
&= \left|(\varphi_{j+1} \circ \cdots \circ \varphi_k)(\eta) - (\varphi_{j+1} \circ \cdots \circ \varphi_k)(z_k)\right| < \frac{\delta_j}{2^k}.
\end{aligned} \tag{3.12}$$

Applying Lemma 2.4, there exists some $0 < \delta_{k+1} < r_{k+1}$ and $x_{k+1} : [0, T_{k+1} - T_k] \times V_{k+1}$, where $V_{k+1} = B(z_{k+1}, \delta_{k+1})$, such that

- for every $\eta \in V_{k+1}$, the function $t \mapsto x_{k+1}(t, \eta)$ is a solution of
$$\dot{x} = -F(T_{k+1} - t, x), \quad x(0) = \eta,$$
for $t \in [0, T_{k+1} - T_k]$;

- the map $\eta \mapsto x_{k+1}(\cdot, \eta)$ is continuous from $V_{k+1}$ into $AC([0, T_{k+1} - T_k], X)$, and in particular, $\varphi_{k+1} : \eta \mapsto x_{k+1}(T_{k+1} - T_k, \eta)$ is continuous on $V_{k+1}$;

- for each $\eta \in V_{k+1}$,
$$|x_{k+1}(t, \eta) - z(T_{k+1} - t)| < \widehat{\delta}_k \leq \delta_k, \qquad \forall t \in [0, T_{k+1} - T_k],$$

and in particular,
$$x_{k+1}(T_{k+1} - T_k, z_{k+1}) \in B(z_k, \widehat{\delta}_k). \tag{3.13}$$

Define inductively, for $0 \leq j \leq k$, $\zeta_{k+1}^{k-j} = \varphi_{k-j+1}(\zeta_{k+1}^{k-j+1})$, where $\zeta_{k+1}^{k+1} = z_{k+1}$. By the choice of $\widehat{\delta}_k$ and (3.13), one sees that
$$\left|\zeta_{k+1}^j - \zeta_k^j\right| = \left|(\varphi_{j+1} \circ \cdots \circ \varphi_k)(\zeta_{k+1}^k) - \zeta_k^j\right| < \frac{\delta_j}{2^k},$$

for all $0 \leq j \leq k-1$. For the case of $j = k$, we have

$$\left|\zeta_{k+1}^k - \zeta_k^k\right| = |x_k(T_{k+1} - T_k, z_{k+1}) - z_k| \leq \widehat{\delta}_k \leq \frac{\delta_k}{2^k}.$$

We have shown that (3.11) holds for all $0 \leq j \leq k$.

Thus, by induction, we conclude that for each $k \geq 1$, there is some $0 < \delta_k < r_k$ and $x_k$ such that

- for each $\eta \in B(z_k, \delta_k)$, $x_k(\cdot, \eta)$ is a solution of

$$\dot{x} \in -F(T_k - t, x), \qquad x(0) = \eta,$$

for $t \in [0, T_k - T_{k-1}]$;

- $x_k$ is continuous from $B(z_k, \delta_k)$ to $AC([0, T_k - T_{k-1}], X)$;

- for each $k \geq 0$, there is a sequence $\{\zeta_j^k\}_{j \geq k}$ such that $\zeta_j^k \in B(z_k, \delta_k) = V_k$ for each $j \geq k$; and

    (a.) $x_{k+1}(T_{k+1} - T_k, \zeta_j^{k+1}) = \zeta_j^k$ for all $j \geq k+1$; and

    (b.) $\left|\zeta_{i+1}^k - \zeta_i^k\right| \leq \frac{\delta_k}{2^i}$ for all $i \geq k$, and hence, $\{\zeta_i^k\}_{i \geq k}$ is a Cauchy sequence.

Since $X$ is a Banach space, $\{\zeta_i^k\}_{i \geq k}$ converges. Let $\zeta^k$ be the limit of the sequence $\{\zeta_i^k\}$. With the continuity property of $x_k$, we get the following compatibility property for the boundary values of $x_k$:

$$x_k(T_k - T_{k-1}, \zeta^k) = \lim_{i \to \infty} x_k(T_k - T_{k-1}, \zeta_i^k) = \lim_{i \to \infty} \zeta_i^{k-1} = \zeta^{k-1}.$$

We now define a trajectory $\gamma(\cdot)$ for the inclusion

$$\dot{x} = F(t, x)$$

by $\gamma(t) = x_k(T_k - t, \zeta^k)$ if $t \in [T_{k-1}, T_k)$. Then $\gamma(0) = \zeta^0 \in B(z_0, r_0)$, and

$$|\gamma(t) - z(t)| \leq r_k$$

for all $t \in [T_k, T_{k+1}]$. This then implies that

$$|\gamma(t) - z(t)| < r(t)$$

for all $t \in [0, T)$. ∎

## 4  An Example

Theorem 3.1 provides a complementary result to the Filippov-Ważewski Theorem in the sense that the approximation of a relaxed trajectory by regular trajectories can be carried out over noncompact intervals. However, as shown in [5] by a counterexample, it is in general impossible to have an infinite-time approximation which satisfies the same initial conditions as a given relaxed trajectory. The counterexample in [5] dealt with a specific function $r(t)$ which decays to 0 as $t \to \infty$. Below we modify the example in [5] to show that even in the case when $r(t) \equiv \varepsilon$ for some $\varepsilon > 0$, it may still be impossible to find an infinite-time approximation with the same initial conditions as a given relaxed trajectories.

**Example 4.1** Consider the system of differential inclusion:

$$\begin{aligned} \dot{x}_1(t) &= x_2^2(t) \\ \dot{x}_2(t) &= x_3^2(t) \\ \dot{x}_3(t) &\in \{-1, 1\}, \end{aligned}$$

and the relaxation to convex values:

$$\begin{aligned} \dot{x}_1(t) &= x_2^2(t) \\ \dot{x}_2(t) &= x_3^2(t) \\ \dot{x}_3(t) &\in [-1, 1]. \end{aligned}$$

Note that $x(t) \equiv 0$ is a relaxed trajectory with initial condition $x(0) = 0$, where $x(t) = (x_1(t), x_2(t), x_3(t))$.

Clearly, the set-valued function $F(x)$ satisfies all conditions in Theorem 3.1. For any $\varepsilon > 0$, applying Theorem 3.1 with $r(t) \equiv \varepsilon$, one sees that there exists a solution $x(t)$ of the original system which satisfies $|x(t)| \leq \varepsilon$ for all $t \geq 0$ with $|x(0)| \leq \varepsilon$.

However, the inclusion does not admit any solution satisfying the condition $|x(t)| \leq \varepsilon$ with $x(0) = 0$. For any solution $x(t)$ with $x(0) = 0$, it holds that

$$x_2(1) = \int_0^1 x_3^2(t)\, dt = \sigma$$

for some $\sigma > 0$. Since $\dot{x}_2(t) \geq 0$ for all $t \geq 0$, one sees that $x_2(t) \geq \sigma$ for all $t \geq 1$. Consequently,

$$x_1(t) = \int_0^t x_2^2(s)\, ds \geq \sigma^2(t-1) \qquad \forall t \geq 1.$$

Hence, it is impossible to have $|x_1(t)| \leq \varepsilon$ for any given $\varepsilon > 0$.  □

## 5  Output Stability

In this section we consider the following system with outputs:

$$\dot{x}(t) \in F(x(t)), \quad y(t) = h(x(t)), \tag{5.14}$$

where $F : X \to \mathcal{P}(X)$ is locally Lipschitz with nonempty compact values, and $h : X \to \mathbb{R}^p$ is continuous. For this system, we will use $y(t)$ to denote $h(x(t))$ when there is no confusion. Throughout this section, we assume that $X$ is finite dimensional.

For each $C \subseteq X$ we let $\mathbf{S}(C)$ denote the set of maximal solutions of (5.14) satisfying $x(0) \in C$. If $C$ is a singleton $\{\xi\}$ we will use the shorthand $\mathbf{S}(\xi)$. We set $\mathbf{S} := \mathbf{S}(X)$, the set of all maximal solutions.

The following notions, which generalize the notions of global asymptotic stability and uniform global asymptotic stability to the output case, were introduced in [4] for systems of differential equations.

**Definition 5.1** The system (5.14) is *output-globally asymptotically stable* (oGAS) if the output trajectory is uniformly stable in the following sense: for all $\varepsilon > 0$, there exists a $\delta > 0$ so that for all $|h(\xi)| \leq \delta$ and all $x(\cdot) \in \mathbf{S}(\xi)$

$$|y(t)| \leq \varepsilon \quad \forall t \geq 0, \tag{5.15}$$

and is output-attractive: for every $x(\cdot) \in \mathbf{S}$,

$$\lim_{t \to \infty} |y(t)| = 0. \tag{5.16}$$

**Definition 5.2** The system (5.14) satisfies the *uniform output-global asymptotic stability condition* (oUGAS) if it is uniformly stable as in (5.15) and the outputs are *uniformly attractive*: for all $\kappa > 0$ and $\varepsilon > 0$, there exists some $T = T(\varepsilon, \kappa) > 0$ so that for all $x(\cdot) \in \mathbf{S}(B(0, \kappa))$,

$$|y(t)| \leq \varepsilon \quad \forall t \geq T. \tag{5.17}$$

The next result generalizes the result [4, Theorem 2] to systems defined by differential inclusions as defined in (5.14).

**Proposition 5.3** Suppose that the system (5.14) is forward complete. Then the system satisfies (oGAS) if and only if it satisfies (oUGAS). □

The proof of Proposition 5.3 follows the same ideas as in the work [4]. A key result used in the proofs is that, for systems that are forward complete, every relaxed trajectory can be approximated by regular trajectories over the interval $[0, \infty)$. With the approximation result for systems of differential inclusions now available, the proofs in [4] trivially apply to Proposition 5.3. To be more specific, an underlying result in the proofs is as in the following.

Given a forward complete system (5.14), a subset $S \subseteq X$, and a maximal solution $x(\cdot)$, we denote the "first crossing time" as

$$\tau(S, x(\cdot)) = \inf\{t \geq 0 : x(t) \in S\}$$

with the convention that $\tau(S, x(\cdot)) = \infty$ if $x(t) \notin S$ for all $t \geq 0$. For a subset $S^o$ of the output space $\mathbb{R}^p$, we denote the "first crossing time" for the output as

$$\tau^0(S, x(\cdot)) = \inf\{t \geq 0 : y(t) \in S^o\}$$

with the convention that $\tau^0(S, x(\cdot)) = \infty$ if $y(t) \notin S^o$ for all $t \geq 0$. Following the same ideas as in the proofs given in Section 2 of [4] with the help Theorem 3.1, one can get the following:

**Lemma 5.4** Let (5.14) be a forward complete system. Assume given a compact subset $C \subseteq X$, an open subset $\Phi \subseteq \mathbb{R}^p$, a compact subset $J \subseteq \Phi$, such that

$$\forall\, x(\cdot) \in \mathbf{S}(C),\ \exists\, t \geq 0 \text{ s.t. } y(t) \in J.$$

Then

$$\sup\{\tau^o(\Phi, x(\cdot)) : x(\cdot) \in \mathbf{S}(C)\} < \infty.$$